\documentclass[12pt,leqno]{article}
\usepackage{amsmath,amsfonts}
\newtheorem{lem}{Lemma}[section]

\newtheorem{theo}[lem]{Theorem}

\newtheorem{prop}[lem]{Proposition}
\newcommand{\vers}{\mathop{\longrightarrow}} 
\newcommand{\R}{{\mathbb R}}
\newcommand{\Sy}{{\mathcal S}}
\newcommand{\tr}{\operatorname{Tr}}

\numberwithin{equation}{section} 

\title{Large deviations for Wishart processes }
\author{C. Donati-Martin\thanks
{Laboratoire de Probabilit\'es et Mod\`eles Al\'eatoires, Universit\'e Paris 6, Site Chevaleret, 13 rue Clisson, F-75013 Paris. email: donati@ccr.jussieu.fr }}
\date{}
\begin{document}
\maketitle
\begin{abstract}
Let $X^{(\delta)}$ be a Wishart process of dimension $\delta$, with values in the set of positive matrices of size $m$. We are interested in the large deviations for a family of matrix-valued processes $\{\delta^{-1} X_t^{(\delta)}, t \leq 1 \}$ as $\delta$ tends to infinity. The process $X^{(\delta)}$ is a solution of a stochastic differential equation with a degenerate diffusion coefficient. Our approach is based upon the introduction of exponential martingales. We give some applications to large deviations for functionals of the Wishart processes, for example the set of eigenvalues.
\end{abstract}


\section{Introduction}
Let $B$ be a $m \times m$ matrix valued Brownian motion. We consider a Wishart process $X_t$, solution of the following SDE, with values in $\Sy_m^{+}$, the set of $m\times m$ real symmetric non-negative matrices:
\begin{equation} \label{sdeW}
dX_t = \sqrt{X_t} \;dB_t + dB'_t \sqrt{X_t} + \delta I_m \;dt, 
\quad X_0= x , 
\end{equation}
where $x \in \Sy_m^{+}$ and $M'$ denotes the transpose of the matrix $M$. \\
We recall the following existence theorem (see M.F. Bru \cite{Bru}): 
\begin{itemize}
\item[] if $\delta \geq m+1$, and $x \in \widetilde{\Sy}_m^{+}$ (the set of positive definite symmetric matrices), then \eqref{sdeW} has a unique strong solution in $\widetilde{\Sy}_m^{+}$.  
\end{itemize}
In fact, we can extend this result to a degenerate initial condition, and in the following, we shall allow $x=0$.

We shall look for a Large Deviation Principle for the $\widetilde{\Sy}_m^{+}$ valued diffusion with small diffusion coefficient:
\begin{equation} \label{eqW}
\left \{ \begin{array}{l}
dX^\epsilon_t = \epsilon ( \sqrt{X^\epsilon_t} \;dB_t + dB'_t \sqrt{X^\epsilon_t} )+ \delta I_m \;dt, \ t \leq T\\
 X^\epsilon_0= x  \end{array} \right.
\end{equation}
with $\delta >0$. For $\epsilon$ small enough, according to the above existence result, \eqref{eqW} has a unique solution $X^\epsilon_t \in \widetilde{\Sy}_m^{+}$ for $t >0$. \\
Note that this problem is equivalent to look for a LDP for the family of processes $\displaystyle (\frac{1}{N} X_t^{(N \delta)}; t \leq 1)$ where $X_t^{(N \delta)}$ denotes a Wishart process of dimension $N \delta$, starting from $Nx$ as $N \vers \infty$. \\
When $m=1$, \eqref{sdeW} is the equation for the squared Bessel process (BESQ) of dimension $\delta$.

In a previous paper \cite{DRYZ}, we studied large deviations for BESQ and squared Ornstein-Uhlenbeck processes. Note that the diffusion coefficient in the BESQ equation is not Lipschitz and the Freidlin-Wentzell theory doesn't apply directly (in the degenerate cases : $x = 0$  or $\delta =0$).We gave three approaches; the first one was based upon exponential martingales, the second one uses the infinite divisibility of the law of BESQ processes (and thus a Cramer theorem) and the third method is a consequence of the continuity of the It\^o map for the Bessel equation (not square), a property proved by Mc Kean \cite{McK}. \\
We also refer to Feng \cite{Fe} for the study of a LDP for squares of Ornstein-Uhlenbeck processes.

\noindent
In the matrix case, due to the restriction on the dimension $\delta$, the laws $Q_x^\delta$ of the Wishart processes are no more infinitely divisible. Moreover, we have no analogue of the Bessel equation for the square root of a Wishart process. \\
Thus, we shall focus on the exponential martingale approach to extend the LDP in the matrix case. Since the delicate point is for a degenerate initial condition, we shall assume that $x = 0$. \\
We denote by $C_0([0,T]; \widetilde{\Sy}_m^{+})$ the space of continuous paths $\varphi_t$ from $[0,T]$ to  $\Sy_m^{+}$ such that $\varphi_0 = 0$ and $\varphi_t \in \widetilde{\Sy}_m^{+}$ for $t>0$.

\noindent
The main result of the paper is:
\begin{theo} \label{theoprinc}
The family $P^\epsilon$ of distributions of $ (X^\epsilon_t ; t \in [0,T])$, solution of \eqref{eqW},  satisfies a LDP in $C_0([0,T]; \widetilde{\Sy}_m^{+})$ with speed $\epsilon^2$  and good rate function 
\begin{equation} \label{taux}
I(\varphi) = \frac{1}{8}  \int_0^T  \tr(k_\varphi (s) \varphi(s) k_\varphi (s)) ds, \quad \varphi \in C_0([0,T]; 
\widetilde{\Sy}_m^{+})
\end{equation}
where $k_\varphi(s)$ is the unique symmetric matrix, solution of
\begin{equation} \label{tauxdef}
 k_\varphi (s) \varphi(s) + \varphi(s) k_\varphi(s)= 2 (\dot{\varphi}(s) - \delta I_m), \, s>0.
 \end{equation}
\end{theo}

\noindent
{\bf Remark:} In the real case ($m=1$), we obtain (see \cite{DRYZ}),
$$ I(\varphi) = \frac{1}{8}  \int_0^T  \frac{(\dot{\varphi}(s) - \delta)^2}{\varphi(s)} ds.$$

The outline of the paper is the following.
In Section 2, we prove an exponential tightness result for the distribution $P^\epsilon$ of $X^\epsilon$.
In section 3, we prove Theorem 1.1 using the approach of exponential martingales. In Section 4, we discuss the Cramer's approach, using the additivity of Wishart processes, when we put some restriction on the parameter $\delta$. In section 5, we give some applications of the contraction principle to obtain 
a  LDP for some functionals of the Wishart process.
\section{Exponential Tightness}
We follow the same lines as in \cite[Section 2]{DRYZ}, that is, we prove exponential tightness in the space $C_\alpha$ of $\alpha$-H\" older continuous functions with $\alpha <1/2$. Let $\alpha <1/2$ and set $\Vert \varphi \Vert_\alpha = \sup_{0 \leq s \not= t \leq T} \frac{\Vert \varphi_t - \varphi_s \Vert}{\vert t- s \vert^\alpha}$ where $\Vert . \Vert$ is a norm on $\Sy_m^{+}$. Since all the norms are equivalent, we shall choose a suitable norm and we consider in this section $
\Vert M \Vert = \sum_{1\leq i,j \leq m} |M_{ij}|$.

\begin{prop}
The family of distributions $P_\epsilon$ of $X^\epsilon$ is exponentially tight in $C_\alpha$, in scale $\epsilon^2$, i.e.
for $L >0$, there exists a compact set $K_L$ in $C_\alpha$ such that:
\begin{equation} \label{tensexpo}
\limsup_{\epsilon \vers 0} \epsilon^2 \ln P( X^\epsilon  \in K_L) \leq -L.
\end{equation}
\end{prop}
{\bf Proof:}
 Let us fix $\alpha' \in (\alpha, 1/2)$ and $R>0$. The closed H\"older ball $B_{\alpha'}(0,R)$
is a compact set   of $C^\alpha([0,1])$. \\ Thus it's enough to estimate $P( \Vert X^\epsilon \Vert_{\alpha'} \geq R)$. For simplicity, we assume $T=1$. 

 $$\Vert X^\epsilon \Vert_{\alpha'} \leq   \Vert M^\epsilon \Vert_{\alpha' } + \delta m$$
 where $M^\epsilon$ is the martingale defined by 
 $$M^\epsilon_t = \epsilon ( \sqrt{X^\epsilon_t} \;dB_t + dB'_t \sqrt{X^\epsilon_t} ).$$

\noindent
{\bf Bounds for $\Vert M^\epsilon \Vert_\alpha$}.
We shall use Garsia-Rodemich-Rumsey's Lemma which asserts that if
$$\int_0 ^1 \int_0 ^1 \Psi \Bigl(\frac{||M^{\epsilon}_t - M^{\epsilon}_s ||}{p(|t-s|)}\Bigr) ds dt \leq K$$
then
$$||M^{\epsilon}_t - M^{\epsilon}_s|| \leq 8 \int_0 ^{\vert t-s \vert} \Psi^{-1} (4K/ u^2 ) dp(u) .$$
Take $\Psi(x) = e ^{c\epsilon^{-2}x} - 1$ for some $0<c<1/2$
and $p(x) = x^{1/2}$. So $\Psi^{-1}(y) = \frac{\epsilon^2}{c} \log (1 +y)$. This yields  (see the same computations in \cite{DRYZ}):
\begin{equation}
P\Big( \| M^{\epsilon} \|_{\alpha'} \geq R\Big)
\leq P\left(\int_0 ^1 \int_0 ^1
\exp \Bigl( c \epsilon^{-2}\frac{||M^{\epsilon}_t - M^{\epsilon}_s ||}{|t-s|^{1/2}}\Bigr) \ ds dt \geq K + 1 \right)  \end{equation}
with $ K = \frac{1}{4} \Bigl( e^{\left(\frac{c\epsilon^{-2}R}{8} - K_2\right) - 4}  - 1 \Bigr) $ and 
$K_2 = 2 \sup _{u \in [0,1]} u^{1/2 - \alpha'} \log \frac{1}{u}$. \\
Now by Markov's inequality,
\begin{equation}\label{MI}
P \Bigl( || M^{\epsilon}||_{\alpha} \geq R \Bigr) \leq 
\frac{1}{K+ 1}\int_0 ^1 \int_0 ^1
E\left[\exp \Bigl(c\epsilon^{-2}\frac{||M^{\epsilon}_t - M^{\epsilon}_s ||}{|t-s|^{1/2}}\Bigr)
\right]ds dt .
\end{equation} 
Now, for a matrix $M$, 
\begin{eqnarray*}
\exp(\lambda ||M||) = \prod_{i,j} \exp(\lambda |M_{ij}|) &\leq &  \prod_{i,j} [\exp(\lambda M_{ij}) + 
\exp(-\lambda M_{ij})] \\
&\leq & m^2 \max [\exp(\lambda M_{ij}) + \exp(-\lambda M_{ij})] \\
&\leq & m^2 \sum_{i,j}  [\exp(\lambda M_{ij}) + \exp(-\lambda M_{ij})]
\end{eqnarray*}
Thus,
\begin{eqnarray*}
\lefteqn{
E[ \exp(\lambda ||M^{\epsilon}_t - M^{\epsilon}_s ||)] } \\
&&\leq  m^2 \sum_{i,j} 
\left(E[ \exp(\lambda (M^{\epsilon}_{i,j}(t) - M^{\epsilon}_{i,j}(s) )] +
E[ \exp(- \lambda (M^{\epsilon}_{i,j}(t) - M^{\epsilon}_{i,j}(s) )] \right)\\
&&\leq 2m^4 \max_{i,j} E[ \exp(2\lambda^2 \langle M^{\epsilon}_{i,j}\rangle_s^t)]
\end{eqnarray*}
where we use in the last inequality the exponential inequality for continuous martingales
$$E[ \exp(\lambda Z_t)] \leq E[\exp(2 \lambda^2\langle Z \rangle_t)].$$
Now,
\begin{eqnarray*}
\langle M^\epsilon_{i,j} \rangle_s^t &= &\epsilon^2 \int_s^t (X_{ii}^\epsilon (u) + X_{jj}^\epsilon (u)) du \\
&\leq & \epsilon^2\int_s^t \tr(X_u^\epsilon) du.
\end{eqnarray*}
Set  $Y^\epsilon_u := \tr(X^\epsilon_u)$, then, $ Y^\epsilon_u$ is a squared Bessel process, solution of the following SDE
\begin{equation} \label{BESQ}
\left\{ \begin{array}{l}
dY^\epsilon_u = 2 \epsilon \sqrt{Y_u^\epsilon} d\beta_u + \delta m \ dt\\
Y_0^\epsilon = 0 
\end{array} \right. \end{equation}
with $\beta$ a real Brownian motion.
Thus, we obtain:
\begin{eqnarray}
E\left[\exp \left( c\epsilon^{-2} \frac{||M^{\epsilon}_t - M^{\epsilon}_s ||}{|t-s|^{1/2}}\right)\right]
&\leq& 2 m^4 \left\{E\left[ \exp
\left(\frac{2c^2 \epsilon^{-2}}{(t-s)}\int_s ^t Y^{\epsilon}_u du\right)\right]\right\}^{1/2} \nonumber\\
&\leq& 2m^4\Bigl\{\frac{1}{t-s} \int_s ^t E\left[\exp \left( 2c^2  \epsilon^{-2}Y^\epsilon_u 
\right)\right] \  du\Bigr\}^{1/2}
\end{eqnarray} 
(by Jensen's  inequality).
Thus, we obtain:
\begin{equation}
\label{PI}
 P(\Vert M^\epsilon \Vert_\alpha \geq R) \leq \frac{2m^4}{K+1}  \left\{ 
 \sup_{u \in [0,1]} E\left[\exp( 2 c^2\epsilon ^{-2} Y^\epsilon_u)\right]  
 \right\}^{1/2}\,,
\end{equation}
 where $K+1 = C \exp(  cR\epsilon^{-2}/8)$ and $C$ a 
 constant. \\
 Now, 
\begin{equation*}
E[\exp( 2 c^2\epsilon ^{-2} Y^\epsilon_u)] =
 Q_{0}^{m \delta \epsilon^{-2}}[\exp( 2 c^2 X_u)]
 \end{equation*}
 where $Q_x^\rho$ denotes the distribution of a squared Bessel process, starting from $x$, of dimension $\rho$. The Laplace transform of the BESQ is known (\cite{RY}) and we obtain: for $c <1/2$,
 $$ Q_{0}^{m \delta \epsilon^{-2}}[\exp( 2 c^2 X_u)] = \left( 1-4c^2 u \right)^
 {- \frac{m \delta\epsilon^{-2}}{2}} .$$
implying
\begin{equation*}
P(\Vert M^\epsilon \Vert_{\alpha'} \geq R) \leq C_m
 A^{ m \delta \epsilon^{-2}} e^{- cR\epsilon^{-2}/8}
\end{equation*}
 for a positive constant $A$. Thus,
 $$\lim_{R \rightarrow +\infty}\limsup_{\epsilon \rightarrow 0} \epsilon^2 
 \ln P(\Vert M^\epsilon \Vert_{\alpha'} \geq R) = - \infty .  \quad \Box
$$

\section{Proof of Theorem \ref{theoprinc}}
From the previous section, we need to prove a weak LDP, that is to prove the upper bound for compact
sets. We assume that $T = 1$. According to \cite{DZ}, we shall prove:
\begin{itemize}
\item[i)] Weak upper bound:
\begin{equation}
\label{weakupper}
\lim_{r \rightarrow 0} \limsup_{\epsilon \rightarrow 0} \epsilon^2 \ln P(X^\epsilon 
\in B_r (\varphi))\leq - I(\varphi) 
\end{equation}
where $B_r (\varphi )$ denotes the open ball with center $\varphi \in {\cal 
C}^{\alpha}_0([0,1]; \widetilde{\Sy}_m^{+})$ and 
radius $r$.
\item[ii)] Lower bound : for all open set $O \subset {\cal 
C}^{\alpha}_0([0,1],; \widetilde{\Sy}_m^{+})$,
\begin{equation}
\label{lower}
 \liminf_{\epsilon \rightarrow 0} \epsilon^2 \ln P(X^\epsilon 
\in O) \geq  - \inf_{\varphi \in O} I(\varphi)\,.
\end{equation}
\end{itemize}

\subsection{The upper bound}
We denote by ${\cal M}_m$, resp. $\Sy_m$ the space of $m\times m$ matrices, resp. symmetric matrices, endowed with the scalar product:
$$ \langle A,B \rangle = \tr(AB').$$
The corresponding norm is denoted by $ \Vert A \Vert_2$.
Set 
$H = \{ h \in  C([0,1]; \Sy_m): \dot h \in L^2([0,1]; \Sy_m) \} $.
For $h \in H$ let
$$M^{\epsilon,h}_t = \exp \left( \frac{1}{\epsilon^2} \{ \int_0^t \tr(h(s) 
(dX^\epsilon_s - \delta I _m \ ds) )\ - \frac{1}{2} \langle Z^\epsilon, Z^\epsilon \rangle_t\} \right), t \leq 1$$
where 
$$Z^\epsilon_t = \int_0^t \tr(h(s) \sqrt{X^\epsilon_s}dB_s +h(s) dB'_s\sqrt{X^\epsilon_s}).$$
$$ \langle Z^\epsilon, Z^\epsilon \rangle_t = 4 \int_0^t \tr(h(s) X^\epsilon_s \ h(s)) ds.$$
$M^{\epsilon,h}$ is a positive, local martingale. In fact, using a Novikov's type criterion (see \cite[Exercise VIII.1.40]{RY}, \cite{DRYZ}), we can prove that $M^{\epsilon,h}$ is a martingale, then, $E(M^{\epsilon,h}_t) =1$.\\
By an integration by parts, we can write:
$$ M^{\epsilon,h}_1 = \exp \left(\frac{1}{\epsilon^2} \Phi(X^\epsilon; h) \right) $$
with
$$ \Phi (\varphi; h)  = G(\varphi; h) - 2 \int_0^1 \tr(h(s) \varphi (s) h(s)) \ ds $$
 and
 $$ G(\varphi ; h) = \tr (h_1 (\varphi_1 - \delta  I_m)) - \int_0^1 \tr((\varphi_s - \delta s I_m) \dot{h}_s) ds$$
for $\varphi \in C_0([0,1]; \widetilde{\Sy}_m^{+})$. \\
{\bf Remark:} If $\varphi$ is absolutely continuous, then, 
$$ G(\varphi; h) = \int_0^1\tr(h(s)  (\dot{\varphi}_s  - \delta I _m \ ds) ). $$
For $\varphi \in C_0([0,1]; \widetilde{\Sy}_m^{+})$, $h \in H$,
\begin{eqnarray*}
P \left(X^\epsilon \in B_r(\varphi)\right) &=& P \left(X^\epsilon \in 
B_r(\varphi); \frac{M^{\epsilon, h}_1}{M^{\epsilon, h}_1}\right) \\
&\leq & \exp \left( - \frac{1}{\epsilon^2} \inf_{\psi \in B_r(\varphi)} 
\Phi(\psi ; h) \right) E(M^{\epsilon,h}_1) \\
&\leq & \exp \left( - \frac{1}{\epsilon^2} \inf_{\psi \in B_r(\varphi)} 
\Phi(\psi ; h) \right)\,,
\end{eqnarray*}
which yields :
$$\limsup_{\epsilon \rightarrow 0} \epsilon^2 \ln P\left(X^\epsilon \in B_r(\varphi)\right) 
\leq - \inf_{ \psi \in B_r(\varphi)} \Phi( \psi; h)\,.$$
For
$h \in H$, the map 
$\varphi \longrightarrow  \Phi(\varphi ;h)$ is continuous on ${\cal C}_0([0,1], \Sy_m^{+})$, so that
$$\lim_{r \rightarrow 0} \limsup_{\epsilon \rightarrow 0} \epsilon^2 \ln P\left(X^\epsilon \in B_r(\varphi)\right) 
\leq -  \Phi(\varphi ; h).$$
Minimizing in $h \in H$, we obtain:
$$\lim_{r \rightarrow 0} \limsup_{\epsilon \rightarrow 0} \epsilon^2 \ln P\left(X^\epsilon \in B_r(\varphi)\right) 
\leq -\sup_{h \in H}  \Phi(\varphi ; h).$$

\begin{prop}
For $\varphi \in C_0([0,1]; \widetilde{\Sy}_m^{+})$,
$$ \sup_{h \in H}  \Phi(\varphi ; h) = I(\varphi) $$
where $I(\varphi)$ is defined by \eqref{taux}.
\end{prop}
{\bf Proof:} Since $\varphi \in C_0([0,1]; \widetilde{\Sy}_m^{+})$,  $\int_0^1 \tr (h_s \varphi_s h_s) ds >0$ for $h \not\equiv 0$.  Replacing $h$ by $\lambda h$ for $\lambda \in \R$, we can see that 
$$ J(\varphi) := \sup_{h \in H} \Phi(\varphi; h) = \frac{1}{8} \sup_{h \in H} \frac{G^2(\varphi;h)}{
\int_0^1 \tr (h_s \varphi_s h_s) ds }.$$
We assume that $J(\varphi) < \infty$.
We denote by $\Vert h \Vert_{L^2(\varphi)}$ the Hilbert norm on $C_0([0,1]; \Sy_m)$ given by 
$$ \Vert h \Vert^2_{L^2(\varphi)} = \int_0^1 \tr (h_s \varphi_s h_s) ds  .$$
The linear form $G_\varphi \ : h \vers G(\varphi; h)$ can be extended to the space $L^2(\varphi)$ and by Riesz theorem, there exists a function $k_\varphi \in L^2(\varphi)$ such that 
\begin{equation}
G(\varphi; h)  = \int_0^1 \tr (h_s \varphi_s k_\varphi(s)) ds 
\end{equation}
Thus, $\varphi$ is absolutely continuous and we have 
\begin{equation} \label{ac}
  \int_0^1 \tr (h_s (\dot{\varphi}_s - \delta I_m)) ds  = \int_0^1 \tr (h_s \varphi_s k_\varphi(s)) ds 
\end{equation}
for all symmetric matrix $h(s)$. Let $k_\varphi$ be given by \eqref{tauxdef}. We refer to the Appendix for the existence of an unique solution of \eqref{tauxdef}. Then, it is easy to see that \eqref{ac} is satisfied for all $h$ symmetric. Moreover, by Cauchy-Schwarz inequality,
\begin{eqnarray*}
\int_0^1 \tr (h_s \varphi_s k_\varphi(s)) ds & \leq & \int_0^1 \tr (h_s \varphi_s h_s)^{1/2} 
 \tr (k_\varphi (s) \varphi_s k_\varphi(s))^{1/2} ds \\
 & \leq &  ( \int_0^1 \tr (h_s \varphi_s h_s) ds )^{1/2}( \int_0^1  \tr (k_\varphi (s) \varphi_s k_\varphi(s)) ds)^{1/2}
 \end{eqnarray*}
 with equality for $h = k_\varphi$. \\
 Thus,
 $\displaystyle \frac{1}{8} \sup_{h \in L^2(\varphi)}  \frac{G^2(\varphi;h)}{\Vert h \Vert^2_{L^2(\varphi)}} = I(\varphi)$. \\
 Now, the equality between $I(\varphi)$ and $J(\varphi)$ follows by density of $H$ in $L^2(\varphi)$.
 $\Box$
\subsection{The lower bound}
In order to prove the lower bound, we first prove 
$$
 \liminf_{\epsilon \rightarrow 0} \epsilon^2 \ln P(X^\epsilon 
\in B_r (\varphi)) \geq - I(\varphi) 
$$
for  all $r >0$ and for $\varphi$ in a subclass ${\cal H}$ of $C_0([0,1]; \widetilde{\Sy}_m^{+})$.
Then, we shall show that this subclass is rich enough.

\noindent
Set ${\cal H}$  the set of functions $\varphi $ such that $I(\varphi) < \infty$ and s.t. $k_\varphi $ defined by \eqref{tauxdef}  belongs to $H$. \\
For $\varphi \in {\cal H}$, set $h_\varphi = \frac{1}{4} k_\varphi$. As in the previous subsection,
we introduce the new probability measure
$$\hat{P} := M^{\epsilon, h_\varphi}_1 P$$ where $P$ is the Wiener measure on $C([0,1]; {\cal M}_{m,m})$.  Under $\hat{P}$,
$$B_t = \hat{B}_t + \frac{2}{\epsilon} \int_0^t (\sqrt{X^\epsilon_s} \ h_\varphi (s) )\ ds$$
where $\hat{B}$ is a Brownian matrix on $\hat{P}$. \\
Thus, under $\hat{P}$, $X^\epsilon$ solves the SDE
$$ dX^\epsilon_t = \epsilon ( \sqrt{X^\epsilon_t} d\hat{B}_t + d\hat{B}'_t  \sqrt{X^\epsilon_t} ) + 
(2( X^\epsilon_t h_\varphi (t) + h_\varphi (t) X^\epsilon_t) + \delta I_m) dt. $$
Under $\hat{P}$, $\displaystyle X^\epsilon_t \vers_{\epsilon \rightarrow 0}^{a.s.} \Psi_t$ solution of 
$$ d\Psi_t = [2(\Psi_t h_\varphi (t) + h_\varphi (t) \Psi_t) + \delta I_m] dt $$
i.e.
$$ \dot{\Psi}_t - \delta I_m = 2(\Psi_t h_\varphi (t) + h_\varphi (t) \Psi_t)  = \frac{1}{2} (  \Psi_t k_\varphi (t) + k_\varphi (t) \Psi_t).$$
Since $k_\varphi$ is continuous, this equation has $\varphi$ as a unique solution; thus 
$$ X^\epsilon_t \vers_{\epsilon \rightarrow 0} \varphi_t \; \hat{P}\  \mbox{a.s.}$$
and
 $\displaystyle \lim_{\epsilon \vers 0} \hat{P}(X^\epsilon \in B_r(\varphi)) = 1$ for every $r>0$. Now,
 \begin{eqnarray*}
P \left(X^\epsilon \in B_r(\varphi)\right) &=& \hat{P} \left(X^\epsilon \in 
B_r(\varphi) \frac{1}{M^{\epsilon, h_\varphi}_1}\right) \\
&\geq & \exp \left( - \frac{1}{\epsilon^2} \sup_{\psi \in B_r(\varphi)} 
F(\psi ; h_\varphi) \right) \hat{P}(X^\epsilon \in B_r(\varphi)) \\
\end{eqnarray*}
which yields :
$$\liminf_{\epsilon \rightarrow 0} \epsilon^2 \ln P\left(X^\epsilon \in B_r(\varphi)\right) 
\geq - \sup_{ \psi \in B_r(\varphi)} F( \psi; h_\varphi)\,$$ and by continuity of $F(.,h)$,
$$ \lim_{r \vers 0} \liminf_{\epsilon \rightarrow 0} \epsilon^2 \ln P\left(X^\epsilon \in B_r(\varphi)\right) 
\geq  F(\varphi; h_\varphi)\, = I(\varphi).$$
We now prove the:
\begin{prop}
For any $\varphi \in C_0([0,1]; \widetilde{\Sy}_m^{+})$ such that $I(\varphi ) < \infty$, there exists a sequence $\varphi_n$ of elements of ${\cal H}$ such that $\varphi_n \vers \varphi$ in $C_0([0,1]; \widetilde{\Sy}_m^{+})$ and $I(\varphi_n) \vers I(\varphi)$.
\end{prop}
{\bf Proof:} We follow the same lines as in the proof of the corresponding result for the scalar case in \cite{DRYZ}.

a) First, let us show that the condition $I(\varphi) < \infty$ implies that 
$$\lim_{t \vers 0} \frac{\varphi_t}{t} = \delta I_m.$$
 From the scalar case, we know that:
\begin{equation} \label{p1}
 \lim_{t \vers 0} \frac{\tr(\varphi_t)}{t} = \delta m.\end{equation}
Indeed, $\tr(X^\varepsilon_t)$ satisfies a LDP (see \eqref{BESQ}) with rate function given by $J( g) =  \frac{1}{8}  \int_0^1 \frac{(\dot{g}(s) - \delta m)^2}{g(s)} ds$ and $J(g) < \infty$ implies that $  \lim_{t \vers 0} \frac{g(t)}{t} = \delta m$. (see \cite{DRYZ}, \cite{Fe}).  From the upper bound, the condition $I(\varphi) < \infty$ implies that  $J(\tr(\varphi)) < \infty$ and thus
\eqref{p1} is satisfied.

\noindent
Let us denote $|| A ||_1 = \tr(|A|)$ and $||A||_2 = (\tr(|A|^2))^{1/2}$ for a matrix $A$.
\begin{eqnarray*}
||\varphi_t - \delta t I_m||_1 &=& \Vert \int_0^t (\dot{\varphi}_s - \delta I_m) ds\,  \Vert_1 \\
&=& \frac{1}{2} \Vert \int_0^t (\varphi_s k_\varphi (s) +  k_\varphi (s)\varphi_s \, \Vert_1 \\
&\leq &  \frac{1}{2}(  \int_0^t \Vert \varphi_s k_\varphi (s)\Vert_1 ds + \int_0^t \Vert k_\varphi (s)\varphi_s \Vert_1 ds ) \\
&\leq &  \int_0^t  \Vert \sqrt{\varphi_s} \Vert_2 \ \Vert \sqrt{\varphi_s} k_\varphi (s)\Vert_2 ds \\
&= &  \int_0^t  (\tr(\varphi_s))^{1/2} (\tr(k_\varphi (s) \varphi_s k_\varphi (s)))^{1/2}  ds \\
&\leq & \left(  \int_0^t  \tr(\varphi_s) ds \right)^{1/2} \left( \int_0^t \tr(k_\varphi (s) \varphi_s k_\varphi (s)) ds \right)^{1/2}  
\end{eqnarray*}
Thus,
$$
||\frac{\varphi_t}{t} - \delta I_m ||_1 \leq \left( \frac{1}{t} \int_0^t \frac{ \tr(\varphi_s)}{s} ds \right)^{1/2}
 \left( \int_0^t \tr(k_\varphi (s) \varphi_s k_\varphi (s)) ds \right)^{1/2}  .$$
 According to \eqref{p1}, the first term in the RHS is bounded and the second tends to 0 as $t$ tends to $0$ since $I(\varphi) < \infty$.
 
 \vspace{.3cm}
 b) As a second step, we approximate $\varphi$ by a function $\psi$ such that $k_\psi \in L^2([0, 1]; {\cal S}_m)$. Set 
 $$ \left\{ \begin{array}{ll}
 \psi_r (t) = \delta t I_m, & t \leq r/2 \\
 \psi_r (t) = \delta r/2\  I_m + (t - r/2) a_r, & r/2 \leq t \leq r \\
 \psi_r (t) = \varphi(t), & t \geq r
 \end{array} \right.
 $$
 where the matrix $a_r$ is chosen such that $\psi$ is continuous in $r$. Let $k_\psi$ the solution of \eqref{tauxdef}  associated with $\psi$. Since $k_\psi(s) =0$ for $s \in [0, r/2]$, and $\psi(s)$ is invertible for $s>0$, $k_\psi \in L^2([0, 1]; {\cal S}_m)$. Obviously, $\displaystyle \psi_r \vers_{r\rightarrow 0} \varphi$ in $C_0([0,1]; \widetilde{\Sy}_m^{+})$.  \\
 It remains to prove the convergence of $I(\psi_r)$ to $I(\varphi)$, or that 
 $$ \int_{r/2}^r \tr(k_\psi (s) \psi(s) k_\psi (s)) ds \vers_{r\rightarrow 0} 0.$$
 \begin{eqnarray*}
 \int_{r/2}^r \tr(k_\psi (s) \psi(s) k_\psi (s)) ds  &= &\int_{r/2}^r \tr(k_\psi (s)(\dot{\psi} (s)- \delta I_m)) ds \\
 &=& \int_{r/2}^r \tr(k_\psi (s)(a_r- \delta I_m)) ds.
 \end{eqnarray*}
 Note that $a_r$ and $k_{\psi_r} (s)$ for $s \in [r/2, r]$ are diagonalisable in the same basis with respective eigenvalues: $(a_i^{(r)})_i$ and $k_i (s) = \frac{a_i^{(r)} - \delta}{\delta r/2 + (s-r/2) a_i^{(r)}}$
 and that, according to step a), $\displaystyle \lim_{r\vers 0} a_i^{(r)} = \delta$.
 Thus, for $r$ small enough,
 $$\int_{r/2}^r \tr(k_\psi (s) \psi(s) k_\psi (s)) ds  = \int_{r/2}^r\sum_i \frac{(a_i^{(r)} - \delta)^2}{
 \delta r/2 + (s-r/2)a_i^{(r)}} ds \leq 1/\delta  \sum_{i=1}^m (a_i^{(r)} - \delta)^2$$
 and the last quantity tends to 0 as $r$ tends to 0.
 
 \vspace{.3cm}
 c) Thanks to b), we must find an approximating sequence $\varphi^{(n)}$ of $\varphi$ in ${\cal H}$ for 
 $\varphi$ satisfying $k_\varphi \in L^2 $. \\
 Let $k^{(n)}$ be a sequence of smooth functions with values in ${\cal S}_m$ such that $k^{(n)}$ converges to $k_\varphi$ in $L^2([0,1], {\cal S}_m)$. Let $\varphi^{(n)}$ be the unique solution of
 $$\left  \{ \begin{array}{l}
 \dot{\varphi}^{(n)}_t  - \delta I_m = k^{(n)}_t \varphi^{(n)}_t + \varphi^{(n)}_t k^{(n)}_t  \\
 \varphi^{(n)}_0 = 0
 \end{array} \right. $$
 Since
 $$||\varphi^{(n)}_t  || \leq \int_0^t \Vert \dot{\varphi}^{(n)}_s \Vert  ds \leq 2  \int_0^t \Vert \varphi^{(n)}_s \Vert \  \Vert k^{(n)}_s \Vert ds
 + \delta,$$
 the Gronwall inequality shows that
 $$ \sup_n \sup_{t \in [0,1]} ||\varphi^{(n)}_t  || < \infty$$ where we have chosen the operator norm on the set of matrices in the previous inequality. Another application of Gronwall's inequality entails that:
 $$ \sup_{t \in[0,1]} || \varphi_t - \varphi^{(n)}_t  || \vers_{n \rightarrow \infty} 0.$$
 Now, the convergence of $I(\varphi^{(n)})$ to $I(\varphi)$ follows from the convergence in $L^2$ of $k^{(n)}$ to $k_\varphi$ and the convergence in $L^\infty([0,1])$ of $\varphi^{(n)}$ to $\varphi$. $\Box$
 \section{The Cramer theorem}
 Let $Q^\delta_x$ denote the distribution on $C(\R, {\cal S}_m^+)$ of the Wishart process of dimension $\delta \geq m+1$, starting from $x \in {\cal S}_m^+$.
 We recall the following additivity property (see \cite{Bru}):
 $$ Q^{\delta}_{x}  \oplus  Q^{ \delta'}_{y} = Q^{\delta+ \delta'}_{x+y}. $$
Let $\delta \geq m+1 $ and take $ \epsilon = \frac{1}{\sqrt{n}}$, then $X^\epsilon$, solution of \eqref{eqW} , is distributed as $\frac{1}{n} \sum_{i=1}^n X_i$ where $X_i$ are independent copies of $Q^\delta_x$. From Cramer's theorem (\cite{DZ} chap. 6), we obtain:
\begin{theo}
Let $\delta \geq m+1$.  The family $P^\epsilon$ of distributions of $ (X^\epsilon_t ; t \in [0,T])$, solution of \eqref{eqW},  satisfies a LDP in $C_0([0,T]; \widetilde{\Sy}_m^{+})$ with speed $\epsilon^2$  and good rate function:
\begin{equation} \label{Fenchel} 
\Lambda^*(\varphi) = \sup_{ \mu \in {\cal M}([0,T], {\cal S}_m)} \left( \int_0^T \tr(\varphi_t d\mu_t) \ - 
\Lambda (\mu) \right),
\end{equation}
where 
\begin{equation} \label{loglaplace} 
\Lambda(\mu) =\ln \left[Q^\delta_x \left( \exp( \int_0^T \tr(X_s d\mu_s))\right)\right].
\end{equation}
\end{theo}
The Laplace transform of the $Q^\delta_x$ distribution can be computed explicitely in terms of Ricatti equation, extending to the matrix case, the well known result for the squared Bessel processes (see \cite{PY}, \cite[Chap. XI]{RY}).
\begin{lem}
Let $\mu$ be a positive ${\cal S}_m^+$-valued measure on $[0,T]$. Then,
\begin{equation} \label{Laplace}
Q^\delta_x \left( \exp(- \frac{1}{2} \int_0^T \tr(X_s d\mu_s))\right) = \exp(\frac{1}{2} \tr(F_\mu (0) x)) 
\exp(\frac{\delta}{2}\int_0^T  \tr(F_\mu (s)) ds) 
\end{equation}
where $F_\mu(s)$ is the ${\cal S}_m$-valued, right continuous  solution of the Riccati equation
\begin{equation} \label{Riccati}
\dot{F}+ F^2 = \mu, \quad F(T) = 0.
\end{equation}
\end{lem}
{\bf Proof:}
From It\^o's formula,
\begin{eqnarray*}
F_\mu(t) X_t &=& F_\mu(0) x + \int_0^t F_\mu(s) dX_s + \int_0^t  dF_\mu(s) X_s\\
&=& F_\mu(0) x + \int_0^t F_\mu(s) dX_s + \int_0^t  d\mu(s) X_s \ - \int_0^t  F_\mu^2(s) X_s\ ds
\end{eqnarray*}
Consider the exponential local martingale
$$Z_t = \exp \left( \frac{1}{2} \int_0^t \tr(F_\mu(s) dM_s) - \frac{1}{2} \int_0^t \tr(F_\mu(s) X_s F_\mu(s)) ds \right)$$
where $M_s = X_s - \delta I_m \ s $. Then,
$$Z_t = \exp \left( \frac{1}{2} \left( \tr(F_\mu(t) X_t) - \tr( F_\mu(0) x) - \delta \int_0^t \tr(F_\mu(s)) ds -
 \int_0^t \tr(X_s d\mu(s))\right)\right).$$
 Now, $X_t$ is positive and $F_\mu(t)$ is negative\footnote{See the Appendix (A.2)}. Thus, $\tr(X_t F_\mu(t))
 \leq 0$ and $Z_t$ is a bounded martingale. The Lemma follows from the equality $E(Z_0) = E(Z_T)$.
 $\Box$
 
 \noindent
 {\bf Remarks:} 
 \begin{enumerate}
 \item The condition $F(T) = 0$ in \eqref{Riccati} is equivalent to $F(T-) = - \mu(\{T\})$.
\item  Taking $d\mu_s = 2 \Theta \delta_1(ds) $ where $\Theta$ is a symmetric positive matrix, we find that
$F_\mu (t) = -2\Theta(I_m + 2(1-t)\Theta)^{-1}, \; t < 1$, from which we recover (see \cite{Bru}):
\begin{equation} \label{LaplaceB}
Q^\delta_x \left( \exp(- \tr(X_1 \Theta))\right) = \det(I_m + 2 \Theta)^{-\delta/2} \exp(-\tr(x(I_m+2 \Theta)^{-1} \Theta)).
\end{equation}
For $m=1$, this example is given in \cite{DRYZ}, Subsection 8.3.
\end{enumerate}

\noindent
Let us try to make the correspondence between $\varphi$ and $\mu$ in \eqref{Fenchel}. If $\mu$ is a negative measure, then, from \eqref{Laplace},
$$\int_0^T \tr(\varphi_t d\mu_t) \ - \Lambda (\mu) = \int_0^T \tr(\varphi_t d\mu_t) \ - 
 \frac{1}{2} \tr(F_{-2\mu} (0) x) - \frac{\delta}{2}\int_0^T  \tr(F_{-2\mu} (s)) ds.$$
 Since $d\mu(t) = - \frac{1}{2} (\dot{F}_t + F^2_t)$, an integration by parts gives:
 \begin{equation} \label{eqF}
 \int_0^T \tr(\varphi_t d\mu_t) \ - \Lambda (\mu)  =  \frac{1}{2}\int_0^T  \tr(F_{-2\mu} (s)(\dot{\varphi}_s - \delta I_m)) ds -  \frac{1}{2}\int_0^T  \tr(F^2_{-2\mu} (s) \varphi_s ) ds.
 \end{equation}
 The optimal function $F(s)$ giving the supremum in \eqref{eqF} solves the equation:
 $$\dot{\varphi}_s - \delta I_m = \varphi_s F(s) + F(s) \varphi_s$$
 that is $F(s) = k_\varphi (s)/2$ where $k_\varphi$ is the solution of \eqref{tauxdef}, and for this $F$, the RHS of \eqref{eqF} is exactly $I(\varphi)$.
 \section{Some applications}
 From the contraction principle, we can obtain a LDP for some continuous functionals of the Wishart process $X^\epsilon$.
 \subsection{The eigenvalues process}
 Let $(\lambda^\epsilon (t) = (\lambda^\epsilon_1 (t), \ldots, \lambda^\epsilon_m (t)); t \in [0,T])$ denote the process of eigenvalues of the process $X^\epsilon$.
 \begin{prop} \label{propev}
 The process $\lambda^\epsilon$ satisfies a LDP in $C_0([0,T], \R_+^m)$ with rate function:
 \begin{equation} \label{tauxev}
 J(x) = \frac{1}{8} \sum_{i=1}^m \int_0^T \frac{(\dot{x}_i (t) - \delta)^2}{x_i(t)} dt.
 \end{equation}
 \end{prop}
 {\bf Remark:} $(\lambda^\epsilon (t))_t$ is solution of the SDE (see \cite{Bru}):
 $$ d\lambda^\epsilon_i (t) = 2 \epsilon \sqrt{\lambda^\epsilon_i (t)} d\beta_i (t) + \{ \delta + \epsilon^2
 \sum_{k\not=i} \frac{ \lambda^\epsilon_i(t) + \lambda^\epsilon_k (t)}{\lambda^\epsilon_i (t) - \lambda^\epsilon_k(t)} \} dt $$
 from which we can guess the form of the rate function $J$ in \eqref{tauxev} since the drift $b_\epsilon$ in the above equation satisfies $\displaystyle b_\epsilon (\lambda) \vers_{\epsilon \rightarrow 0} \delta$. Nevertheless, since the drift
 $b_\epsilon (\lambda)$ explodes on the hyperplanes $\{ \lambda_i = \lambda_j\}$ and the diffusion coefficient is degenerate, the classical results (see \cite[Theo. V.3.1]{FW}) do not apply.
 
 \vspace{.3cm}
 \noindent
 {\bf Proof of the Proposition \ref{propev}}
 According to the contraction principle,
 $$ J(x) = \inf\{ I(\varphi); e.v.(\varphi) = x \}.$$
 Write $\varphi_t = P_t^{-1} \Lambda_t P_t$ where $\Lambda$ is the diagonal matrix of eigenvalues of $\varphi_t$ and $P_t$ is an orthogonal matrix. Then,
 $$\dot{\varphi}_t =  P_t^{-1} \dot{\Lambda}_t P_t +  \dot{P_t^{-1}} \Lambda_t P_t + 
  P_t^{-1} \Lambda_t \dot{P}_t.$$
  We denote by $\tilde{k}_t $ the matrix $P_t k_\varphi(t) P_t^{-1}$ where $k_\varphi$ solves \eqref{tauxdef}. Then,
  $$\tr(k_\varphi (t) \varphi(t) k_\varphi(t)) = \tr(\tilde{k}_t \Lambda(t) \tilde{k}_t)$$
  and
  $$\tilde{k}_{ij} (t) \lambda_i(t) + \tilde{k}_{ij} (t) \lambda_j(t) = 2(\dot{\lambda}_i(t) - \delta) \delta_{ij} +
  R_{ij}(t)$$
  where the matrix $R$ is defined by 
  $$ R(t) = P_t  \dot{P_t^{-1}} \Lambda_t +  \Lambda_t \dot{P}_t P^{-1}_t.$$
  Now, it is easy to verify that $R_{ii} (t) =0$, thus:
  $$\tr(\tilde{k}_t \Lambda(t) \tilde{k}_t) = \sum_{i} \frac{(\dot{\lambda}_i (t) - \delta)^2}{\lambda_i(t)} +
  \sum_{i \not= j} \frac{R_{ij}^2(t) \lambda_j(t)}{\lambda_i(t) + \lambda_j(t)}$$
  and the infimum of the above quantity is obtained for $R\equiv 0$, corresponding to $P_t$ independent of $t$. For this choice, $I(\varphi) = J(\lambda)$ where $\lambda$ is the set of e.v. of $\varphi$. $\Box$
 \subsection{A LDP for the r.v. $X^\epsilon_1$}
 \begin{prop}
 The r.v. $X^\epsilon_1$ satisfies a LDP, in scale $\epsilon^2$, with rate function:
 \begin{equation} \label{tauxrv}
 K(M) =  \frac{1}{2} \tr(M) -  \frac{\delta}{2}\ln(\det(M)) -  \frac{m \delta}{2} +  \frac{m\delta}{2} \ln(\delta), \; M \in {\cal S}_m^{+}.
 \end{equation}
 \end{prop}
 {\bf Remark:} For $m =1$, 
 $$K(a) = \frac{1}{2} [ (a - \delta) - \ln(a/\delta)], \ a >0$$ 
 which corresponds (for $\delta =1$) to the rate function obtained in the study of a LDP for a $\chi_2(n)$ distribution as $n \vers \infty$.
 
  \vspace{.3cm}
 \noindent
 {\bf Sketch of proof:} 
 
 i) Since the application $\varphi \vers \varphi(1)$ is continuous, we must minimize $I(\varphi)$ under the constraint $\varphi(1) = M$.
 The optimal path $\varphi$   solves the  Euler Lagrange equation (see \cite{L}, Chap. 7), given in terms of $k_\varphi$ by:
 $$ 2 \dot{k}_\varphi (s) + k_\varphi^2(s) = 0, s \in (0,1).$$
 This leads to  $k_\varphi^{-1} (t) = \frac{t}{2} I_m + C$ and $\varphi(t) = \delta t I_m + t^2 A$ with a matrix $A$ determined by $\varphi(1) = M$. Note that this is the same path as in Section 4, Remark 2. \\
 Now, it is easy to verify that for $\varphi(t) = \delta t I_m +
 t^2 (M - \delta I_m)$, $I(\varphi) = K(M)$ where $K$ is given by \eqref{tauxrv}.
 
 ii) Of course, we can compute $K$ directly, using the Laplace transform \eqref{LaplaceB} (with $x=0$) and then,
 $$K(M) = \sup_{\Theta} \{\tr(\Theta M) + \frac{\delta}{2} \ln(\det(I_m - 2 \Theta))\}.$$
 The optimal $\Theta_0$ is given by $M = \delta(I_m - 2\Theta_0)^{-1}$. $\Box$
\subsection{ A LDP for the largest eigenvalue}
Let us denote by $\lambda^\epsilon_{max}$ the largest eigenvalue of the Wishart process $X^\epsilon$.
\begin{prop}
The process $\{ \lambda^\epsilon_{max}(t), t \in [0,T] \}$ satisfies a LDP in $C_0([0,T; \R_+)$ with rate function given by
$$I_{max} (f) = \inf \{ J(x), x = (f, x_2, \ldots, x_m), x_i(t) \leq f(t) \mbox{ for } i= 2, \ldots m \}$$
where $J$ is given by \eqref{tauxev}. \\
For $f$ belonging to a class of functions ${\cal F}$ to be defined in the proof,
\begin{equation} \label{tauxmax}
 I_{max}(f) = \frac{1}{8} \left[ \int_0^T \frac{(\dot{f}_t - \delta)^2}{f_t} dt + (m-1) \int_0^T 
\frac{(\dot{\underline{f}}_t - \delta)^2}{\underline{f}_t} \right]
\end{equation}
where $\underline{f}(t) = \delta t + \inf_{s \leq t} (f(s) - \delta s)$.
\end{prop}

\noindent
{\bf Proof:} 
We assume that the eigenvalues are given in decreasing order: $\lambda_1(t) \geq \lambda_2(t) \geq \ldots \geq \lambda_m (t)$. \\
According to the contraction principle, $I_{max}$ is given by the minimium of :
$$J(x) =  \frac{1}{8} \sum_{i=1}^m \int_0^T \frac{(\dot{x}_i (t) - \delta)^2}{x_i(t)} dt$$
under the constraint $\{ x_i(t) \leq f(t), i = 2, \ldots m\}$ with $x_1 = f$ fixed. \\
Set 
$$ F(y) =  \frac{1}{8} \int_0^T \frac{(\dot{y} (t) - \delta)^2}{y(t)} dt;$$
 $F$ is a convex function on $C_0([0,T); \R_+)$ and introduce the convex function $G_f(y) = y -f \in C([0,T); \R)$. \\
The problem is to minimize $F(y)$ under the constraint $G(y) \leq 0$.

\noindent
To $f$, we associate the measure $\mu_f$ associated to the Ricatti equation
$$2 \mu_f = \dot{H} + H^2 \mbox{ on } (0,T), \ H(T) = - 2\mu_f(T)$$
with $\displaystyle H_t =  \frac{(\dot{f} (t) - \delta)}{2f(t)} $. \\
Then, we define the measure $d\tilde{\mu}_f (t) = 
d\mu_f (t) 1_{(\underline{f} (t) = f(t))}$. \\

Let ${\cal F} = \{ f; d\tilde{\mu}_f  \mbox{ is a positive measure on $[0,T]$} \}$. For $f \in {\cal F}$, let us show that  the Lagrangian 
$$ L(y, \mu) = F(y) +\langle G_f(y), \mu \rangle$$
has a saddle point at $(\underline{f}, \tilde{\mu}_f)$, i.e.,
\begin{equation} \label{saddle}
L(\underline{f}, \mu) \leq L(\underline{f}, \tilde{\mu}_f) \leq L(y, \tilde{\mu}_f).
\end{equation}
for all $y  \in C_0([0,T); \R_+)$ and all positive measure $\mu$.\\
The first inequality follows from
$$\langle G_f(\underline{f}), \mu \rangle \leq 0 = \langle G_f(\underline{f}), \tilde{\mu}_f \rangle $$
since $supp(\tilde{\mu}_f ) \subset \{t, f(t) = \underline{f}(t) \}$. \\
For the second inequality, we must show that $\underline{f}$ minimize $F(y) + \langle G_f(y), \tilde{\mu}_f \rangle$. 
The optimal path of this problem of minimization  solves the Euler- Lagrange equation (see \cite{DRYZ}):
\begin{equation} \label{EL}
\frac{d}{dt} \left( \frac{\partial g}{\partial b}(y, \dot{y}) \right) =  \frac{\partial g}{\partial a}(y, \dot{y})
+ \tilde{\mu}_f  \mbox{ on }  (0,T), \ \left( \frac{\partial g}{\partial b}(y, \dot{y}) \right)_{t = T} = - \tilde{\mu}_f (T).
\end{equation}
with $g(a,b) = \frac{(b- \delta)^2}{8a}$.
The auxiliary function $H_t = \frac{(\dot{y} (t) - \delta)}{2y(t)} $ associated to the optimal path $y$ satisfies the Ricatti equation:
$$2 \tilde{\mu}_f = \dot{H} + H^2; \; H(T) = - 2\tilde{\mu}_f(T).$$
By the choice of $\tilde{\mu}_f$, it is easy to see that $\underline{f}$ solves the Euler-Lagrange equation
\eqref{EL} (or the associated Ricatti equation). 

\noindent
According to Luenberger (Theo2,  Section 8.4), the existence of this saddle point implies that :
$$ \underline{f} \mbox{ minimize }F(y) \mbox{ under the constraint } G_f(y) \leq 0 . \; \Box$$ 

\noindent
For a fixed time, we have the following result:
\begin{prop}
The r.v. $ \lambda^\epsilon_{max}(1)$ satisfies a LDP in $ \R_+$ with rate function given by
\begin{equation} \label{tauxmax1>}
 K_{max}(a) =\frac{a}{2} - \frac{\delta}{2} \ln(a) - \frac{\delta}{2} + \frac{\delta}{2} \ln(\delta) \mbox{ if } a > \delta
 \end{equation}
 \begin{equation} \label{tauxmax1<}
 K_{max} (a) =m\left( \frac{a}{2} - \frac{\delta}{2} \ln(a) - \frac{\delta}{2} + \frac{\delta}{2} \ln(\delta) \right) \mbox{ if } a \leq \delta
 \end{equation}
\end{prop}
The proof is immediate from \eqref{tauxrv}. We minimize $K(M)$ under the constraint $||M|| = a$, where $||.||$ denotes the operator norm.

\section{Appendix}
{\bf (A.1) On the equation AX+XA = B.}\\

\noindent
Let $A$ and $B$ two symmetric matrices, with $A$ strictly positive. We are looking for a symmetric matrix $X$, solution of the equation
( see \eqref{tauxdef}):
$$ AX +XA = B \quad (*)$$
Since $A$ is symmetric, let $P$ and $D$ be  orthogonal and positive diagonal matrices such that $A = P^{-1}DP$. Then, from $(*)$, the symmetric matrix $\tilde{X} = PXP^{-1}$ satisfies:
$$ D\tilde{X} + \tilde{X} D = PBP^{-1}:= \tilde{B}$$ 
that is:
$$ d_i \tilde{X}_{ij} + \tilde{X}_{ij} d_j = \tilde{B}_{ij}$$
and $\displaystyle \tilde{X}_{ij} = \frac{\tilde{B}_{ij}}{d_i + d_j}$. Thus, $X$ is uniquely determined.

\vspace{.5cm}
\noindent
{\bf (A.2) On the Riccati equation.} \\

\noindent
We consider the Ricatti equation (see \eqref{Riccati}):
\begin{equation} \label{Riccatibis}
\dot{F}+ F^2 = \mu, \quad F(T) = 0.
\end{equation}
or
$$F(t) = C + \mu(]0,t]) - \int_0^t F^2(s) ds$$
where $C$ is chosen that $F(T)= 0$.
We diagonalize $F(t)$: $F_t = P^{-1}_t D_t P_t$  with $D_t$ the matrix of eigenvalues of $F_t$ and $P_t$ orthogonal. Then, the Ricatti equation can be written as:
$$\dot{D} (t) + D^2(t) = P(t) \mu_t P^{-1}(t) + R_t$$
where $R$ is a matrix, whose diagonal entries are zeroes. Set $\nu = P \mu P^{-1}$, then $\nu$ is a positive ${\cal S}_m^+$-valued measure and the eigenvalues of $F$ satisfy the scalar Riccati equation:
$$ \dot{d_i}(t)+ d_i^2(t) = \nu_{ii}(t), \; d_i(T) = 0$$
with $\nu_{ii}$ a positive measure on [0,T]. We know (see \cite{RY}, Chapter XI) that $d_i(t) \leq 0$ ($d_i$ is related to the decreasing solution of the Sturm Liouville equation $\phi''_i = \phi_i \nu_{ii}$).
It follows that the matrix $F(t)$ is symmetric negative.



\begin{thebibliography}{99}
\bibitem{Bru}  Bru, M.-F.: Wishart processes. {\it J}. {\it Theo}. {\it Probab}., {\bf 4} (1991), 725--751.
\bibitem{DZ} Dembo, A. and Zeitouni, O.: {\it Large deviations techniques and applications}. Second Edition, Springer, 1998.
\bibitem{DRYZ} Donati-Martin, C.; Rouault, A.; Yor, M. and Zani, M.: Large deviations for squares of Bessel and Ornstein-Uhlenbeck processes. {\it Prob. Th. Rel. Fields} {\bf 129} (2004) 261-289.
\bibitem{Fe} Feng, S.:  The behaviour near the boundary of some degenerate diffusions under random perturbations. In {\it Stochastic models} (Ottawa, ON, 1998), 115-123, Providence (2000) Amer. Math. Soc.
\bibitem{FW} Freidlin, M.I. and Wentzell A.D.: {\it Random Perturbations of Dynamical Systems}.  Springer-Verlag, New York, 1984.
\bibitem{L} Luenberger, D.G.: {\it Optimization by vector space methods.}  John Wiley, 1969.
\bibitem{McK} Mc Kean, H. P. : The Bessel motion and a singular integral equation. {\it Mem. Coll. Sci. Univ. 
Kyoto. Ser. A, Math.} {\bf 33} (1960) 317-322.
\bibitem{PY} Pitman, J.  and Yor, M.:  A decomposition of Bessel bridges. {\it Z. W} {\bf 59} (1982) 425-457.
\bibitem{RY} Revuz, D.  and  Yor, M.: {\it Continuous martingales and Brownian motion}, {\it 3rd Ed}., 
Springer, Berlin, 1999.
\end{thebibliography}
\end{document}